\documentclass[12pt,a4paper]{amsart}

\usepackage[latin1]{inputenc}
\usepackage{amsfonts,amsthm,amssymb,latexsym,amsmath}
\usepackage[all,cmtip]{xy}
\usepackage{multirow}
\usepackage{enumerate}
\usepackage[left=3cm,right=2cm,top=3cm,bottom=2cm]{geometry}\usepackage[square,sort,comma,numbers]{natbib}

\newtheorem{Teo}{Theorem}

\newtheorem{Lema}[Teo]{Lemma}

\usepackage[
	pdfauthor={Charles}
	pdftitle={Report I}
	pdfdisplaydoctitle=true,
	pdfstartview={FitH},
	bookmarks=false,
]{hyperref}

\title{Semistable rank 2 sheaves with singularities of mixed
dimension on $\mathbb{P}^3$}

\author[A. N. Ivanov]{Alexey N. Ivanov}
\address{Department of Mathematics\\
National Research University
Higher School of Economics\\
6 Usacheva Street\\
119048 Moscow, Russia}
\email{anivanov$_{-}$1@edu.hse.ru}

\author[A. Tikhomirov]{Alexander S. Tikhomirov}
\address{Department of Mathematics\\
National Research University
Higher School of Economics\\
6 Usacheva Street\\
119048 Moscow, Russia}
\email{astikhomirov@mail.ru}

\begin{document}

\begin{abstract}
\noindent We describe new irreducible components of the Gieseker-Maruyama moduli scheme $\mathcal{M}(3)$ of semistable rank 2 coherent sheaves with Chern classes $c_1=0, c_2=3, c_3=0$ on $\mathbb{P}^3$, general points of which correspond to sheaves whose singular loci contain components of dimensions both 0 and 1. These sheaves are produced by elementary transformations of stable reflexive rank 2 sheaves with $c_1=0,\ c_2=2$ along a disjoint union of a projective line and a collection of points in $\mathbb{P}^3$. The constructed  families of sheaves provide first examples of irreducible components of the Gieseker-Maruyama moduli scheme such that their general sheaves have singularities of mixed dimension.


\end{abstract}

\maketitle

\section{Introduction}\label{intro}

\vspace{5mm}

Let $\mathcal{M}(0,k,2n)$ be the Gieseker-Maruyama moduli scheme of semistable rank-2 sheaves with Chern classes $c_1=0,\ c_2=k,\ c_3=2n$ on the projective space $\mathbb{P}^3$. Denote $\mathcal{M}(k)=\mathcal{M}(0,k,0)$.
By the singular locus of a given
$\mathcal{O}_{\mathbb{P}^3}$-sheaf $E$ we understand
the set $\mathrm{Sing}(E)=\{x\in\mathbb{P}^3\ |\ E$ is not locally free at the point $x\}$. $\mathrm{Sing}(E)$ is always a proper closed subset of $\mathbb{P}^3$
and, moreover, if $E$ is a semistable sheaf of nonzero
rank, every irreducible component of $\mathrm{Sing}(E)$ has dimension at most 1.

In this article we study the Gieseker-Maruyama moduli scheme $\mathcal{M}(3)$. In \cite{ES}, \cite{JMT2} it was shown that the scheme $\mathcal{M}(3)$ contains at least the following irreducible components:
\begin{enumerate} [(i)]
\item  the instanton component $\mathcal{I}$ of dimension 21 whose general point
\footnote{Here and below by a general point of a given irreducible set we understand a closed point belonging to its open dense subset. }
corresponds to a locally free instanton sheaf \cite{JMT2};

\item the Ein component $\mathcal{N}$ of dimension 21 whose general point corresponds to a locally free sheaf that is the cohomology sheaf of a monad of the form
$$ 0 \longrightarrow \mathcal{O}_{\mathbb{P}^3}(-2) \longrightarrow \mathcal{O}_{\mathbb{P}^3}(-1) \oplus 2 \mathcal{O}_{\mathbb{P}^3} \oplus \mathcal{O}_{\mathbb{P}^3}(1) \longrightarrow \mathcal{O}_{\mathbb{P}^3}(2) \longrightarrow 0;$$

\item the component $\mathcal{C}$ of dimension 21 whose general point corresponds to a sheaf $E$ having singularities along a smooth plane cubic $C$ and satisfying the following exact triple
$$ 0 \longrightarrow E \longrightarrow 2 \mathcal{O}_{\mathbb{P}^3} \longrightarrow L(2) \longrightarrow 0,$$
\noindent where $L$ is a line bundle over $C$ with Hilbert polynomial $P_L(k)=3k$ such that $L^{\otimes 2} \not\simeq \omega_C$;

\item the four components $\mathcal{T}(s),~s=1,...,4,$ of dimensions $21+4s$, respectively, whose general points correspond to sheaves $E$ having singularities along the sets of $s$ distinct points $\{q_1,...,q_s\}\subset \mathbb{P}^3$ and satisfying the exact triple
$$ 0 \longrightarrow E \longrightarrow F \longrightarrow \bigoplus\limits_{i=1}^{s} \mathcal{O}_{q_i}  \longrightarrow 0,$$
\noindent where $F$ is a reflexive stable rank 2 sheaf with Chern classes $c_1=0,\ c_2=3,\ c_3=2s$ such that $\text{Sing} (F) \cap \{q_1,...,q_s\} = \emptyset$ (in \cite{JMT2} these components are denoted by $\mathcal{T}(3,s)$).
\end{enumerate}
Note that general points of these components correspond to sheaves with singularities of pure dimension, i. e. such that all irreducible components of singular loci of these sheaves have the same dimensions.
 The main result of the present article, Theorem \ref{Thm1} below, gives
a description of new irreducible components of $\mathcal{M}(3)$ whose general points correspond to sheaves with singularities of mixed dimension, i. e. sheaves whose singular loci contain components of dimensions both 0 and 1. For this purpose, we define 7 families of sheaves $\mathcal{X}(n,s)$ from  $\mathcal{M}(3)$ and prove that, among their closures  $\overline{\mathcal{X}(n,s)}$ in $\mathcal{M}(3)$, there are 3 irreducible components of $\mathcal{M}(3)$, namely, $\overline{\mathcal{X}(1,0)},\  \overline{\mathcal{X}(2,0)}$, and $\overline{\mathcal{X}(2,1)}$ (see Theorem \ref{Thm1}). Moreover, in Section \ref{section4} we prove that closures of the rest 4 families $\overline{\mathcal{X}(n,s)}$ are proper subsets of the known components $\mathcal{T}(s),~ s=1,2,3,$ of $\mathcal{M}(3)$, so they do not constitute components of $\mathcal{M}(3)$ (see Theorem \ref{Thm2}).

Remark that, up to now, there were known only those components of the Gieseker-Maruyama schemes  $\mathcal{M}(k)$, general points of which are semistable sheaves with singularities of pure dimension. Namely, in \cite{JMT2} there were constructed the first infinite series of components of the schemes $\mathcal{M}(k)$, with $k$ growing, such that their general points are sheaves with singularities of pure dimension 0, respectively, 1. Thus, Theorem \ref{Thm1} provides first examples of irreducible components of the moduli scheme $\mathcal{M}(k)$ for $k=3$, the general sheaves of which have singularities of mixed dimension. 

The families $\mathcal{X}(n,s)$ are constructed in the following way.
Let $\mathcal{R}(0,2,2n)$ be the moduli scheme of stable reflexive coherent rank 2 sheaves on $\mathbb{P}^3$ with Chern classes $c_1=0,\ c_2=2,\ c_3= 2n$, where $n=1,\ 2$.  $\mathcal{R}(0,2,2n)$ is an open subset of the  Gieseker-Maruyama moduli scheme $\mathcal{M}(0,2,2n)$, see \cite{H-r}.
Let $\mathrm{Sym}^s(\mathbb{P}^3)^*$ be a dense open subset of
$\mathrm{Sym}^s(\mathbb{P}^3)$ consisting of disjoint unions of $s$
distint points in $\mathbb{P}^3$. We consider any set of $s$ distinct points
$$
W=\{q_1,...,q_s\}\in\mathrm{Sym}^s(\mathbb{P}^3)^*,\ \ \ 0\le s\le n+1,
$$
as a reduced scheme. (Here for $s=0$ we set by definition $W=\emptyset$.)

For any point $[F] \in \mathcal{R}(0,2,2n)$, by Grauert-M{\"u}llich Theorem (see \cite{HL}), for $0\le s \le n + 1$, the set
$$
X_{[F]}=\{(l,W)\in\text{Gr}(2,4)\times \mathrm{Sym}^s(\mathbb{P}^3)^*\ |\ l \cap W = \emptyset, (l\sqcup W)\cap\mathrm{Sing}(F)=\emptyset\ \mathrm{and}\
F|_l \simeq 2 \mathcal{O}_l \}
$$
is a dense open subset of $\text{Gr}(2,4)\times \mathrm{Sym}^s(\mathbb{P}^3)^*$.
For $u=(l,W)\in X_{[F]}$ set
$$
r=n+1-s,\ \ \ Q_{u,r}:=\mathcal{O}_l(r) \oplus \mathcal{O}_W,
$$
and consider the open dense subset $\text{Hom}(F, Q_{u,r})_e \subset \text{Hom}(F, Q_{u,r})$ consisting of epimorphisms $\phi:F\twoheadrightarrow Q_{u,r}$. For any element $\phi \in \text{Hom}(F, Q_{u,r})_e$ we can take the kernel of $\phi$:
$$
E=\ker~\phi.
$$
(This procedure of passing from $F$ to $E$ is sometimes called an \textit{elementary transformation of $F$ along $l\sqcup W$}.)
It is easy to see that $E$ is a stable sheaf and defines the point $[E]$ of the scheme $\mathcal{M}(3)$ (use the condition $r+s=n+1,\ r,s \geq 0$). Besides, $\ker \phi \simeq \ker \phi'$ if and only if there is an automorphism $g \in \text{Aut}(Q_{u,r})$ such that $\phi'=g \circ \phi$.
Denote by $[\phi]$ the equivalence class of $\phi$ modulo
$\text{Aut}(Q_{u,r})$.
Now consider the sets
$$
\widetilde{\mathcal{X}}(n,s)=\{x=([F],u,[\phi_x])~|~ [F] \in \mathcal{R}(0,2,2n),\ u \in X_{[F]},\ [\phi_x] \in \text{Hom}(F, Q_{u,r})_e / \text{Aut}(Q_{u,r}) \}.
$$
Since $\mathcal{R}(0,2,2n)$ are reduced irreducible schemes (see \cite{Chang}),
the sets $\widetilde{\mathcal{X}}(n,s)$ have a natural structure of reduced irreducible schemes. Furthermore, by the above, there is a well-defined injective morphism
$$
f:\widetilde{\mathcal{X}}(n,s) \longrightarrow \mathcal{M}(3),~~ x \mapsto [\ker~\phi_x].
$$
\noindent
Set
$$
\mathcal{X}(n,s):=f(\widetilde{\mathcal{X}}(n,s))
\subset\mathcal{M}(3),\ \ \ r+s=n+1,\ \ \  r,s \geq 0.
$$
The dimensions of the varieties $\mathcal{X}(n,s)$ are calculated by the following formula
$$
\text{dim~}\mathcal{X}(n,s)=\text{dim~}\widetilde{\mathcal{X}}(n,s)=\text{dim}~\mathcal{R}(0,2,2n)+\text{dim~} \text{Gr}(2,4) + \text{dim~}\mathrm{Sym}^s(\mathbb{P}^3) +
$$
$$
+\text{dim~}\text{Hom}(F, Q_{u,r})_e / \text{Aut}(Q_{u,r})=13+4+ 3 s + 2 (r+1+s) - (1+s)=18+2r+4s.
$$

A sufficient condition for the variety $\overline{\mathcal{X}(n,s)}$ to be an irreducible component of $\mathcal{M}(3)$ is the equality
$$
\text{dim~Ext}^1(E,E)=\text{dim~}\mathcal{X}(n,s) \text{~for any point~} [E] \in \mathcal{X}(n,s).
$$
\noindent In Section \ref{section3}\ it will be shown that this equality is satisfied for families $\mathcal{X}(1,0),\ \mathcal{X}(2,0)$ and $\mathcal{X}(2,1)$, this giving the proof of Theorem \ref{Thm1}. From the construction of the families it is clear that general points of these irreducible components correspond to sheaves with singularities along $l \sqcup \text{Sing}(F) \sqcup W$. Since $\text{dim Sing}(F)=0$, these singularities have mixed dimension.

Throughout this work, the base field  $\mathbf{k}$ is an algebraically closed field of characteristic 0. Also, for simplicity, we will not
distinguish between a stable rank-2 $\mathcal{O}_{\mathbb{P}^3}$-sheaf $E$
and its isomorphism class $[E]$ as a point in the Gieseker-Maruyama moduli scheme.

{\bf Acknowledgements.} ANI was supported in part by the Simons Foundation. AST was supported by a subsidy to the HSE from the Government of the Russian Federation for the implementation of Global Competitiveness Program. AST also acknowledges the support from the Max Planck Institute for Mathematics in Bonn, where this work was finished during the winter of 2017.

\vspace{5mm}

\section{Some properties of stable reflexive rank 2 sheaves with $c_1=0, c_2=2$ on $\mathbb{P}^3$}\label{section 2}

\vspace{5mm}

According to \cite[p. 63, 66]{Chang}, a general sheaf $[F]\in\mathcal{R}(0,2,2n)$, where $n=1$ or 2, satisfies the exact triple
\begin{equation}\label{1eqn16}
\xi:\ \ \ 0 \longrightarrow \mathcal{O}_{\mathbb{P}^3}(-1)\longrightarrow F\longrightarrow\mathcal{I}_{Y}(1)\longrightarrow 0,
\end{equation}
where for $n=1$ the scheme $Y=m\sqcup C_2$ is a disjoint union of a reduced line $m$ and a smooth conic $C_2$ in $\mathbb{P}^3$; respectively,  for $n=2$ the scheme $Y=C_3$ is a smooth twisted cubic in $\mathbb{P}^3$.
Moreover, the extension $\xi\in \text{Ext}^1(\mathcal{I}_{Y}(1),\mathcal{O}_{\mathbb{P}^3}(-1))$ under isomorphisms
$$
\text{Ext}^1(\mathcal{I}_{Y}(1),\mathcal{O}_{\mathbb{P}^3}(-1))\simeq H^0(\mathcal{E}xt^1(\mathcal{I}_{Y}(1),
\mathcal{O}_{\mathbb{P}^3}(-1)))\simeq H^0(\mathcal{E}xt^2(\mathcal{O}_{Y}(2),\mathcal{O}_{\mathbb{P}^3}))$$
corresponds to the global section $\sigma \in H^0(\mathcal{E}xt^2(\mathcal{O}_{Y}(2),\mathcal{O}_{\mathbb{P}^3}))$ such that $(\sigma)_0=\{x_1,...,x_{2n}\}$ is a union of $2n$ distinct points of the component of $Y$ which is not a line, and
\begin{equation}\label{1eqn18}
\mathcal{E}xt^1(F,\mathcal{O}_{\mathbb{P}^3})\simeq \bigoplus\limits_{i=1}^{2n}\mathcal{O}_{x_i}.
\end{equation}

\begin{Lema}\label{lemma 1}
The sheaf $\mathcal{H}om(F,F)$ satisfies the following exact triple
\begin{equation*}\label{1eqn25}
0\longrightarrow F(-1)\longrightarrow \mathcal{H}om(F,F)\longrightarrow\mathcal{G} \longrightarrow0,
\end{equation*}
where the sheaf $\mathcal{G}$ fits in an exact triple
\begin{equation*}\label{triple for G}
0\to F\otimes\mathcal{I}_{Y}(1)\to\mathcal{G}\to Q\to0,
\ \ \ \ \ \ \dim Q=0.
\end{equation*}
\end{Lema}
\noindent\textit{Proof.} Applying the functor $\mathcal{H}om(-,F)$ to (\ref{1eqn16}) we obtain the exact sequence
\begin{equation}\label{1eqn31}
0\longrightarrow\mathcal{H}om(\mathcal{I}_{Y}(1),F)\longrightarrow\mathcal{H}om(F,F)\longrightarrow F(1) \overset{\delta}{\longrightarrow}
\end{equation}
$$
\overset{\delta}{\longrightarrow}\mathcal{E}xt^1(\mathcal{I}_{Y}(1),F)\longrightarrow\mathcal{E}xt^1(F,F).
$$
Since $\mathrm{hd}(F)=1$, there is a locally free $\mathcal{O}_{\mathbb{P}^3}$-resolution
\begin{equation}\label{1eqn32}
0\longrightarrow L_1\longrightarrow L_0\longrightarrow F\longrightarrow0
\end{equation}
of the sheaf $F$. Applying the functor $\mathcal{H}om(\mathcal{I}_{Y}(1),-)$ to (\ref{1eqn32}) and considering an easily verifiable equality $\mathcal{E}xt^2(\mathcal{I}_{Y}(1),L_1)=0$, we obtain an exact triple
\begin{equation}\label{1eqn33}
\mathcal{E}xt^1(\mathcal{I}_{Y}(1),L_1)\longrightarrow\mathcal{E}xt^1(\mathcal{I}_{Y}(1),L_0)\longrightarrow
\mathcal{E}xt^1(\mathcal{I}_{Y}(1),F)\longrightarrow 0.
\end{equation}
An easy computation shows that the sheaves $\mathcal{E}xt^1(\mathcal{I}_{Y}(1),L_0)$ and $\mathcal{E}xt^1(\mathcal{I}_{Y}(1),L_1)$ are locally free $\mathcal{O}_{Y}$-sheaves. Hence, $\mathcal{E}xt^1(\mathcal{I}_{Y}(1),F)$ is $\mathcal{O}_{Y}$-sheaf and, moreover, it is
generically a rank 2 locally free $\mathcal{O}_{Y}$-sheaf. Furthermore, the morphism $\delta$ in (\ref{1eqn31}) factors through the morphism of restriction $\otimes\mathcal{O}_{Y}:F(1)\longrightarrow F(1)\otimes\mathcal{O}_{Y}$. Since $F(1)\otimes\mathcal{O}_{Y}$ is generically a locally free $\mathcal{O}_{Y}$-sheaf of rank 2, this yields  $\delta=j\circ(\otimes\mathcal{O}_{Y})$, where $j$ is a generically injective morphism of $\mathcal{O}_{Y}$-sheaves. Consequently, since
$\ker(\otimes\mathcal{O}_{Y})=F\otimes\mathcal{I}_{Y}(1)$, it follows that the sheaf $\mathcal{G}:=\ker\delta$
satisfies the exact triple
\begin{equation}\label{1eqn34}
0\longrightarrow\mathcal{H}om(\mathcal{I}_{Y}(1),F)\longrightarrow\mathcal{H}om(F,F)\longrightarrow \mathcal{G}\longrightarrow0
\end{equation}
and the triple (\ref{triple for G}).

Now show that
\begin{equation}\label{1eqn35}
\mathcal{H}om(\mathcal{I}_{Y}(1),F)\simeq F(-1).
\end{equation}
Indeed, applying the functor $\mathcal{H}om(-,F)$ to
\begin{equation}\label{1eqn19}
0\longrightarrow\mathcal{I}_{Y}\longrightarrow\mathcal{O}_{\mathbb{P}^3}\longrightarrow\mathcal{O}_{Y}\longrightarrow0
\end{equation}
we obtain an exact sequence
\begin{equation}\label{1eqn36}
0\longrightarrow F(-1)\longrightarrow\mathcal{H}om(\mathcal{I}_{Y}(1),F)\longrightarrow\mathcal{E}xt^1(\mathcal{O}_{Y}(1),F).
\end{equation}
On the other hand, applying the functor $\mathcal{H}om(\mathcal{O}_{Y}(1),-)$ to (\ref{1eqn32}) yields an exact triple
$$
...\longrightarrow\mathcal{E}xt^1(\mathcal{O}_{Y}(1),L_0)\longrightarrow\mathcal{E}xt^1(\mathcal{O}_{Y}(1),F)\longrightarrow
\mathcal{E}xt^2(\mathcal{O}_{Y}(1),L_1)\longrightarrow...
$$
Since $F$ is locally free outside the points $\{x_1,...,x_{2n}\}$, the sheaf $\mathcal{E}xt^1(\mathcal{O}_{Y}(1),F)$ either equals to zero or is an artinian sheaf.
On the other hand,  $\mathcal{E}xt^1(\mathcal{O}_{Y}(1),L_0)=0$ and the sheaf $\mathcal{E}xt^2(\mathcal{O}_{Y}(1),L_1)$ is locally free $\mathcal{O}_{Y}$-sheaf, since $L_0$ and $L_1$ are locally free $\mathcal{O}_{\mathbb{P}^3}$-sheaves. Hence, the sheaf $\mathcal{E}xt^1(\mathcal{O}_{Y}(1),F)$ is a subsheaf of $\mathcal{E}xt^2(\mathcal{O}_{Y}(1),L_1)$, so it cannot be a non-zero artinian sheaf. Consequently, $\mathcal{E}xt^1(\mathcal{O}_{Y}(1),F)=0$ and (\ref{1eqn35}) follows from (\ref{1eqn36}).\\
The statement of Lemma 1 now follows from (\ref{1eqn34}) and (\ref{1eqn35}).                            \hfill$\Box$

\begin{Lema}\label{lemma 2}
Let $[F] \in \mathcal{R}(0,2,2n)$, i. e. $F$ is a stable reflexive rank 2 sheaf on $\mathbb{P}^3$ with Chern classes $c_1=0,\ c_2=2,\ c_3=2n$, where $n=1,2$. Then the following equality holds
\begin{equation*}\label{1eqn43}
\ \ \ h^2(\mathcal{H}om(F,F))=0.
\end{equation*}
\end{Lema}
\noindent\textit{Proof.} Applying the functor $-\otimes\mathcal{O}_{Y}(1)$ to (\ref{1eqn16}) we obtain the exact sequence
\begin{equation}\label{1eqn39}
\mathcal{T}or_1(\mathcal{I}_{Y}(2),\mathcal{O}_{Y})\overset{\varepsilon_1}{\longrightarrow}\mathcal{O}_{Y}\longrightarrow F(1)\otimes\mathcal{O}_{Y}\overset{\varepsilon_0}{\longrightarrow}\mathcal{I}_{Y}(2)\otimes\mathcal{O}_{Y}\longrightarrow0.
\end{equation}
It is easy to see the sheaves $F(1)\otimes\mathcal{O}_{Y}$ and $\mathcal{I}_{Y}(2)\otimes\mathcal{O}_{Y}$ are rank-2 $\mathcal{O}_{Y}$-sheaves and $\varepsilon_0$ is an epimorphism. So $\ker\varepsilon_0$ is artinian sheaf and the morphism $\varepsilon_1$ is non-zero on each component of the curve $Y$. Next, standard computation gives the following
$$\mathcal{T}or_1(\mathcal{I}_{Y}(2),\mathcal{O}_{Y})\simeq\mathcal{T}or_2(\mathcal{O}_{Y}(2),\mathcal{O}_{Y})\simeq\det N_{{Y}/\mathbb{P}^3}^{\vee}(2)
\simeq \left\{
\begin{array}{cl}
\mathcal{O}_m\oplus\mathcal{O}_{C_{2}}(-1) & \text{~for~} n=1, \\
\mathcal{O}_{C_{3}}(- 4 pt) & \text{~for~} n=2. \\
\end{array}
\right.$$
\noindent So we obtain $\mathrm{coker~}\varepsilon_1=\mathrm{ker~}\varepsilon_0\simeq \bigoplus\limits_{i=1}^{2n} \mathcal{O}_{x_i}$ due to $\mathrm{Sing}(F)=\{x_1,...,x_{2n}\}$ that is the following triple is exact
$0\longrightarrow \bigoplus\limits_{i=1}^{2n} \mathcal{O}_{x_i} \longrightarrow F(1)\otimes\mathcal{O}_{Y}\longrightarrow\mathcal{I}_{Y}(2)\otimes\mathcal{O}_{Y}\longrightarrow0$.
This triple as a triple of the $\mathcal{O}_{Y}$-sheaves on the smooth curve $Y$ is splitted. Thus we have the following isomorphism
\begin{equation}\label{1eqn38}
F(1)\otimes\mathcal{O}_{Y}\simeq \bigoplus\limits_{i=1}^{2n} \mathcal{O}_{x_i} \oplus\mathcal{I}_{Y}(2)\otimes\mathcal{O}_{Y}.
\end{equation}

Next, we have
\begin{equation}\label{1eqn30}
\mathcal{I}_{Y}(2)\otimes\mathcal{O}_{Y}\simeq N_{{Y}/\mathbb{P}^3}^{\vee}(2)\simeq \left\{
\begin{array}{cl}
2\mathcal{O}_m(1)\oplus\mathcal{O}_{C_{2}}\oplus\mathcal{O}_{C_{2}}(1) & \text{~for~} n=1, \\
\mathcal{O}_{C_{3}}(pt) \oplus \mathcal{O}_{C_{3}}(pt) & \text{~for~} n=2. \\
\end{array}
\right.
\end{equation}

\noindent So (\ref{1eqn38}) and (\ref{1eqn30}) yield
\begin{equation}\label{1eqn40}
\ \ \ h^i(F(1)\otimes\mathcal{O}_{Y})=0,\ \ \ i\ge1.
\end{equation}

\noindent It is easy to see that $\mathcal{T}or_1 (F, \mathcal{O}_{Y}(1))=0$ because $\mathrm{hd}(F)=1$, so the sheaf $F\otimes\mathcal{I}_{Y}$ is torsion free and the following triple is exact
\begin{equation}\label{1eqn24}
0\longrightarrow F\otimes\mathcal{I}_{Y}(1)\longrightarrow F(1)\longrightarrow F\otimes\mathcal{O}_{Y}(1)\longrightarrow0.
\end{equation}

Besides, according to \cite[Table 2.8.1, 2.12.2]{Chang} we have
\begin{equation}\label{1eqn42}
h^1(F(1))=h^2(F(1))=h^2(F(-1))=h^3(F(-1))=0
\end{equation}
From (\ref{1eqn40})-(\ref{1eqn42}) and Lemma \ref{lemma 1} it follows that
$$
\ \ \ h^2(\mathcal{H}om(F,F))=h^2(\mathcal{G})=h^2(F(1) \otimes \mathcal{I}_{Y})=h^1(F(1) \otimes \mathcal{O}_{Y})=0.
$$
\hfill$\Box$

\vspace{5mm}

\section{Components of $\mathcal{M}(3)$ formed by sheaves with singularities of mixed dimension} \label{section3}

\vspace{5mm}

\begin{Teo}\label{Thm1}
The families $\overline{\mathcal{X}(1,0)}, \ \overline{\mathcal{X}(2,0)},\ \overline{\mathcal{X}(2,1)}$ are irreducible components of the moduli scheme $\mathcal{M}(3)$, of dimensions 22, 24, 26, respectively.
\end{Teo}

\vspace{1mm}\noindent\textit{Proof.} Consider an arbitrary sheaf $[E] \in \mathcal{X}(n,s)$ which according to the construction satisfies the following exact triple
\begin{equation}\label{202eqn48}
0 \longrightarrow E \overset{\varphi}{\longrightarrow} F \longrightarrow \mathcal{O}_{l}(r) \oplus \mathcal{O}_W \longrightarrow 0,
\end{equation}
where $l$ is a line in $\mathbb{P}^3$, $[F] \in \mathcal{R}(0, 2, 2n)$, $W=\{q_1,...,q_s\}$ is a reduced subscheme of $s$ distinct points in $\mathbb{P}^3$ such that
\begin{equation}\label{202eqn49}
l \cap W = \emptyset, \ \ \ \mathrm{Sing}(F) \cap (l \sqcup W) = \emptyset, \ \ \ F|_l \simeq 2 \mathcal{O}_l.
\end{equation}
Since $\mathrm{hd}(F)=1$ and $\dim\mathrm{Sing}(F)=0$, we have
\begin{equation}\label{202eqn47}
\mathcal{E}xt^{\ge 1}(\mathcal{O}_{\mathbb{P}^3},F)
=\mathcal{E}xt^{\ge 2}(F,F)=0,\ \ \ \dim\mathcal{E}xt^1(F,\mathcal{O}_{\mathbb{P}^3})=\dim\mathcal{E}xt^1(F,F)=0.
\end{equation}
From (\ref{202eqn49}), (\ref{202eqn47}) it follows that
\begin{equation}\label{202eqn51}
\mathcal{H}om(F,\mathcal{O}_{l}(r)\oplus \mathcal{O}_W)\simeq2\mathcal{O}_{l}(r)\oplus 2 \mathcal{O}_W
\end{equation}
and the morphism $\varphi$ in (\ref{202eqn48}) induces the isomorphism of the artinian sheaves
\begin{equation}\label{202eqn52}
{}^1\varphi:\ \mathcal{E}xt^1(F,E)\overset{\simeq}{\to}\mathcal{E}xt^1(F,F).
\end{equation}
For the same reason
\begin{equation}\label{202eqn53}
\mathcal{H}om(\mathcal{O}_{l}(r)\oplus \mathcal{O}_W,F)=\mathcal{E}xt^1(\mathcal{O}_{l}(r)\oplus \mathcal{O}_W,F)=0,
\end{equation}
\begin{equation}\label{202eqn6a}
\mathcal{H}om(F,\mathcal{O}_{l}(r)\oplus \mathcal{O}_W)\simeq2\mathcal{O}_{l}(r)\oplus 2 \mathcal{O}_W,\ \ \ \mathcal{E}xt^i(F,\mathcal{O}_{l}(r) \oplus \mathcal{O}_W)=0,\ \ \ i\ge1.
\end{equation}
Similarly we have the isomorphism of the artinian sheaves
\begin{equation}\label{202eqn54}
\mathcal{E}xt^1(F,F)\simeq\mathcal{E}xt^1(F|_U,F|_U)\overset{\simeq}{\longrightarrow}\mathcal{E}xt^1(E|_U,F|_U).
\end{equation}
\noindent due to the isomorphism
$$
E|_U\simeq F|_U,\ \ \ U:=\mathbb{P}^3\setminus (l \cup W),
$$
\noindent that follows from (\ref{202eqn48}).

The inclusion of the artinian sheaf $\mathcal{E}xt^1(F,F)$ into the sheaf $\mathcal{E}xt^1(E,F)$
\begin{equation}\label{202eqn55}
\mathcal{E}xt^1(F,F)\hookrightarrow\mathcal{E}xt^1(E,F)
\end{equation}
\noindent follows from (\ref{202eqn54}) as a direct summand such that
\begin{equation}\label{202eqn56}
\mathcal{E}xt^1(E,F)\simeq\mathcal{E}xt^1(F,F)\oplus\mathcal{E}xt^1(E,F)|_{U'},\ \ \ U'=\mathbb{P}^3\setminus\mathrm{Sing}(F).
\end{equation}
Note that the sheaf $F|_{U'}$ is locally free, so
\begin{equation}\label{202eqn56a}
\mathcal{E}xt^1(E,F)|_{U'}\simeq\mathcal{E}xt^1(E|_{U'},F|_{U'})\simeq\mathcal{E}xt^1(E|_{U'},\mathcal{O}_{U'})\otimes F|_{U'}
\end{equation}
Besides, $\mathcal{E}xt^2(\mathcal{O}_{l}(r)\oplus \mathcal{O}_W,\mathcal{O}_{U'})\simeq(\det N_{l/\mathbb{P}^3})(-r)\simeq\mathcal{O}_l(2-r)$. Applying the functor $\mathcal{H}om(-,\mathcal{O}_{U'})$ to the exact triple (\ref{202eqn48}) restricted on $U'$ and taking into account that the sheaves $\mathcal{E}xt^1(F|_{U'},\mathcal{O}_{U'})$ and $\mathcal{E}xt^2(F|_{U'},\mathcal{O}_{U'})$ are vanished because the sheaf $F|_{U'}$ is locally free we have $\mathcal{E}xt^1(E|_{U'},\mathcal{O}_{U'})\simeq\mathcal{O}_l(2-r)$. So $\mathcal{E}xt^1(E|_{U'},\mathcal{O}_{U'})\otimes F|_{U'}\simeq2\mathcal{O}_l(2-r)$ due to (\ref{202eqn49}). The isomorphism $\mathcal{E}xt^1(E,F)|_{U'}\simeq2\mathcal{O}_l(2-r)$ follows from that and (\ref{202eqn56a}). Hence, from (\ref{202eqn56}) we have the isomorphism
\begin{equation}\label{202eqn56b}
\mathcal{E}xt^1(E,F)\simeq\mathcal{E}xt^1(F,F)\oplus2\mathcal{O}_l(2-r).
\end{equation}

Next, in view of (\ref{202eqn53}) applying the functor $\mathcal{H}om(\mathcal{O}_{l}(r)\oplus \mathcal{O}_W,-)$ to the triple(\ref{202eqn48}) yields the isomorphism
\begin{equation}\label{202eqn57}
\mathcal{E}xt^1(\mathcal{O}_{l}(r)\oplus \mathcal{O}_W,E)\simeq\mathcal{H}om(\mathcal{O}_{l}(r)\oplus \mathcal{O}_W,\mathcal{O}_{l}(r)\oplus \mathcal{O}_W)\simeq\mathcal{O}_{l}\oplus \mathcal{O}_W.
\end{equation}
Besides, applying the functor $\mathcal{H}om(F,-)$ to (\ref{202eqn48}) we have
\begin{equation}\label{202eqn58}
0 \longrightarrow \mathcal{H}om(F,E) \longrightarrow \mathcal{H}om(F,F) \longrightarrow \mathcal{H}om(F,\mathcal{O}_{l}(r)\oplus \mathcal{O}_W) \overset{\partial}{\longrightarrow}\mathcal{E}xt^1(F,E)\overset{{}^1\varphi}{\longrightarrow}\mathcal{E}xt^1(F,F)
\end{equation}
From here in view of (\ref{202eqn51}) and (\ref{202eqn52}) it follows that the following triple is exact
\begin{equation}\label{202eqn59}
0 \longrightarrow \mathcal{H}om(F,E) \longrightarrow \mathcal{H}om(F,F) \longrightarrow2\mathcal{O}_{l}(r)\oplus 2 \mathcal{O}_W\longrightarrow0.
\end{equation}
Next, applying the functor $\mathcal{H}om(-,E)$ to (\ref{202eqn48}) and taking into account (\ref{202eqn55}) and (\ref{202eqn57}) we have the following triple
\begin{equation}\label{202eqn60}
0 \longrightarrow \mathcal{H}om(F,E) \longrightarrow \mathcal{H}om(E,E) \longrightarrow\mathcal{O}_{l}\oplus \mathcal{O}_W\longrightarrow0.
\end{equation}
The monomorphism $E\hookrightarrow F$ in (\ref{202eqn48}) induces the monomorphism $\tau:\mathcal{H}om(E,E) \longrightarrow \mathcal{H}om(F,F)$ which along with the triples (\ref{202eqn59}) and (\ref{202eqn60}) is included into the following commutative diagram
\begin{equation}\label{diag}
\xymatrix{
&  &  & 0 \ar[d] \\
  & 0 \ar[d] & 0 \ar[r] \ar[d] & \mathcal{O}_{l}\oplus \mathcal{O}_W \ar[r]^-{=} \ar[d] & \\
0 \ar[r] & \mathcal{H}om(F,E) \ar[r] \ar[d] & \mathcal{H}om(F,F) \ar[r] \ar[d] & 2\mathcal{O}_{l}(r)\oplus 2 \mathcal{O}_W \ar[r] \ar[d] & 0 & \\
0 \ar[r] & \mathcal{H}om(E,E) \ar[r]{\tau} \ar[d] & \mathcal{H}om(E,F) \ar[r] \ar[d] & \mathrm{coker}~\tau \ar[r] \ar[d] & 0 &\\
\ar[r]^-{=} & \mathcal{O}_{l} \oplus \mathcal{O}_W \ar[r] \ar[d] & 0  & 0  & &\\
  & 0 &  &}
\end{equation}

\noindent
It is easy to see that $\mathrm{coker}~\tau |_{l}$ is a locally free $\mathcal{O}_{l}$-sheaf (see \cite[proof of Lemma 5.1]{JMT1}), so the right vertical triple in this diagram yields
\begin{equation}\label{202eqn61}
\mathrm{coker}~\tau\simeq\mathcal{O}_{l}(2r)\oplus \mathcal{O}_W.
\end{equation}
and the following triple is exact
\begin{equation}\label{202eqn62}
0 \longrightarrow \mathcal{H}om(E,E) \longrightarrow \mathcal{H}om(E,F) \longrightarrow \mathcal{O}_{l}(2r) \oplus \mathcal{O}_W \longrightarrow 0.
\end{equation}
Applying the functor $\mathcal{H}om(E,-)$ to (\ref{202eqn48}) and taking into account (\ref{202eqn62}) we obtain the exact sequence
\begin{equation}\label{202eqn62a}
\begin{aligned}
 0 \longrightarrow\mathcal{O}_l(2r) \oplus \mathcal{O}_W & \longrightarrow \mathcal{H}om(E,\mathcal{O}_l(r) \oplus \mathcal{O}_W)\longrightarrow\mathcal{E}xt^1(E,E) \longrightarrow \mathcal{E}xt^1(E,F) \longrightarrow \\ & \longrightarrow \mathcal{E}xt^1(E, \mathcal{O}_l(r) \oplus \mathcal{O}_W) \overset{\psi}{\longrightarrow} \mathcal{E}xt^2(E,E)
 \end{aligned}
\end{equation}

\noindent On the other hand, applying the functor $\mathcal{H}om(-,\mathcal{O}_l(r)\oplus \mathcal{O}_W)$ to (\ref{202eqn48}) and taking into account (\ref{202eqn6a}) and the standard isomorphisms
$$\mathcal{H}om(\mathcal{O}_l(r)\oplus \mathcal{O}_W,\mathcal{O}_l(r)\oplus \mathcal{O}_W)\simeq\mathcal{O}_l\oplus \mathcal{O}_W,$$
$$\mathcal{E}xt^1(\mathcal{O}_l(r)\oplus \mathcal{O}_W,\mathcal{O}_l(r)\oplus \mathcal{O}_W)\simeq N_{l/\mathbb{P}^3} \oplus 3 \mathcal{O}_W\simeq2\mathcal{O}_l(1) \oplus 3 \mathcal{O}_W,$$
$$\mathcal{E}xt^2(\mathcal{O}_l(r)\oplus \mathcal{O}_W,\mathcal{O}_l(r)\oplus \mathcal{O}_W)\simeq\det N_{l/\mathbb{P}^3} \oplus 3 \mathcal{O}_W\simeq\mathcal{O}_l(2) \oplus 3 \mathcal{O}_W,$$
we obtain the exact sequence
\begin{equation}\label{202eqn62b}
\begin{aligned}
0\longrightarrow\mathcal{O}_l \oplus \mathcal{O}_W \longrightarrow2\mathcal{O}_l(r) \oplus 2 \mathcal{O}_W\overset{\gamma}{\longrightarrow} \mathcal{H}om(E,\mathcal{O}_l(r) \oplus \mathcal{O}_W) \longrightarrow2\mathcal{O}_l(1) \oplus 3 \mathcal{O}_W\longrightarrow0
 \end{aligned}
\end{equation}
and the isomorphism
\begin{equation}\label{202eqn62c}
\mathcal{E}xt^1(E,\mathcal{O}_l(r)\oplus \mathcal{O}_W)\simeq\mathcal{E}xt^2(\mathcal{O}_l(r)\oplus \mathcal{O}_W,\mathcal{O}_l(r)\oplus \mathcal{O}_W)\simeq\mathcal{O}_l(2) \oplus 3 \mathcal{O}_W.
\end{equation}
It is easy to see that $\mathrm{im}\gamma |_{l}$ is locally free $\mathcal{O}_l$-sheaf in (\ref{202eqn62b}), so $\mathrm{im}\gamma$ is isomorphic to $\mathcal{O}_l(2r) \oplus \mathcal{O}_W$. Hence, the following triple is exact
\begin{equation}\label{202eqn62h}
0\longrightarrow\mathcal{O}_l(2r) \oplus \mathcal{O}_W \longrightarrow\mathcal{H}om(E,\mathcal{O}_l(r) \oplus \mathcal{O}_W)\longrightarrow2\mathcal{O}_l(1) \oplus 3 \mathcal{O}_W \longrightarrow 0.
\end{equation}

\noindent Besides, applying the functor $\mathcal{H}om(-,E)$ to (\ref{202eqn48}) we obtain the isomorphism
\begin{equation}\label{202eqn62j}
\mathcal{E}xt^2(E,E) \simeq \mathcal{E}xt^3(\mathcal{O}_l(r) \oplus \mathcal{O}_W,E) \simeq \mathcal{E}xt^3(\mathcal{O}_W,E),
\end{equation}
\noindent since, in view of (\ref{202eqn47}), $\mathcal{E}xt^{\ge 2}(F,E) \simeq \mathcal{E}xt^{\ge 2}(F,F) \simeq 0$ and $\text{hd}(\mathcal{O}_l)=2$. On the other hand, from (\ref{202eqn56b}) it follows that $\mathcal{E}xt^1(E,F)|_{W} = 0$. So $\ker \psi|_{W} = 0$ and, in view of (\ref{202eqn62j}), (\ref{202eqn62c}), we have the isomorphism
\begin{equation}\label{202eqn62k}
\ker \psi \simeq \mathcal{O}_l(2).
\end{equation}

Thus from (\ref{202eqn62a})-(\ref{202eqn62k}) follows the exact sequence
\begin{equation}
\begin{aligned}\label{202eqn62d}
0\longrightarrow 2\mathcal{O}_l(1) \oplus 3 \mathcal{O}_W\longrightarrow \mathcal{E}xt^1(E, E) \longrightarrow\mathcal{E}xt^1(F,F)\oplus2\mathcal{O}_l(2-r) \longrightarrow \mathcal{O}_l(2) \longrightarrow 0.
\end{aligned}
\end{equation}

\noindent Since $\mathcal{E}xt^1(F,F)$ is an artinian sheaf, the following triple is exact
$$
0\longrightarrow2\mathcal{O}_l(1) \oplus 3 \mathcal{O}_W \longrightarrow\mathcal{E}xt^1(E,E)\longrightarrow\mathcal{E}xt^1(F,F)\oplus\mathcal{O}_l(2-2r)
\longrightarrow0.
$$
From this triple we obtain
\begin{equation}\label{202eqn62e}
\begin{aligned}
h^0(\mathcal{E}xt^1(E,E)) & =\left\{
\begin{array}{cl}
4+3s+h^0(\mathcal{E}xt^1(F,F)) & \text{~for~} r \geq 2, \\
7+3s-2r+h^0(\mathcal{E}xt^1(F,F)) & \text{~for~} r = 0 \text{~or~} 1.
\end{array}\right.\ \
\end{aligned}
\end{equation}

Next, the triple (\ref{202eqn62}) in view of the diagram (\ref{diag}) can be written by the following way
\begin{equation}\label{202eqn90}
0\longrightarrow\mathcal{H}om(E,E)\longrightarrow\mathcal{H}om(F,F)\longrightarrow\mathcal{O}_l(2r) \oplus \mathcal{O}_W\longrightarrow0.
\end{equation}
Since $E$ and $F$ are stable, they are simple, i. e. $h^0(\mathcal{H}om(E,E))=h^0(\mathcal{H}om(E,E))=1$. Therefore from Lemma \ref{lemma 2} and (\ref{202eqn90}) it follows that
\begin{equation}\label{202eqn91}
h^1(\mathcal{H}om(E,E))=1+2r+s+h^1(\mathcal{H}om(F,F)),\ \ \ h^2(\mathcal{H}om(E,E))=0.
\end{equation}

Besides, from Lemma \ref{lemma 2} and \cite[Theorem 2.8, 2.12]{Chang} we have
\begin{equation}\label{202eqn92}
h^0(\mathcal{E}xt^1(F,F))+h^1(\mathcal{H}om(F,F))=\text{dim Ext}^1(F,F)=\text{dim~}\mathcal{R}(0,2,2n)=13.
\end{equation}
Consequently, from (\ref{202eqn62e}), (\ref{202eqn91}), (\ref{202eqn92}) and the exact sequence
$$
0\longrightarrow H^1(\mathcal{H}om(E,E))\longrightarrow\mathrm{Ext}^1(E,E)\longrightarrow H^0(\mathcal{E}xt^1(E,E))\longrightarrow  H^2(\mathcal{H}om(E,E))
$$
we obtain
$$
\dim\mathrm{Ext}^1(E,E)=\left\{
\begin{array}{cl}
18 + 4s + 2r & \text{~for~} r \geq 2, \\
21 + 4s & \text{~for~} r = 0 \text{~or~} 1.
\end{array}\right.
$$

\noindent Thus we have the following table of dimensions from which the statement of the theorem follows.

\begin{center}
\begin{tabular}{|c|c|c|}
\hline
The family & Dimension of the family & $\text{dim~Ext}^1(E,E)$  \\
\hline
$\mathcal{X}(1,0)$ & 22  & 22  \\
\hline
$\mathcal{X}(1,1)$ & 24  & 25  \\
\hline
$\mathcal{X}(1,2)$ & 26  & 29  \\
\hline
$\mathcal{X}(2,0)$ & 24  & 24  \\
\hline
$\mathcal{X}(2,1)$ & 26  & 26  \\
\hline
$\mathcal{X}(2,2)$ & 28  & 29  \\
\hline
$\mathcal{X}(2,3)$ & 30  & 33  \\
\hline
\end{tabular}
\end{center}

Since $\mathrm{Ext}^1(E,E)$ is the Zariski tangent space to $\mathcal{M}(3)$ at the point $[E]$, the coincidence of dimensions of the irreducible families $\mathcal{X}(1,0)$,
$\mathcal{X}(2,0)$ and $\mathcal{X}(2,1)$ with the dimensions of the Zariski tangent spaces to $\mathcal{M}(3)$
at their general points yields Theorem \ref{Thm1}.
\hfill$\Box$

The rest families $\mathcal{X}(1,1)$, $\mathcal{X}(1,2),$ $\mathcal{X}(2,2)$, $\mathcal{X}(2,3)$ in the above table will be considered in the next Section.

\section{Deformations of sheaves} \label{section4}

In this Section we will consider the families $\mathcal{X}(1,1)$, $\mathcal{X}(1,2),$ $\mathcal{X}(2,2)$, $\mathcal{X}(2,3)$ and will show that their closures are proper subsets of the known components $\mathcal{T}(s), ~ s=1,2,3,$ of $\mathcal{M}(3)$ - see Theorem \ref{Thm2}.

As it was mentioned in Section \ref{intro}, a general sheaf $E$ of the component $\mathcal{T}(s)$ of $\mathcal{M}(3),\  s=1,2,3$ is defined by the following triple

\begin{equation}\label{E in T1}
0 \longrightarrow E \longrightarrow F \longrightarrow \mathcal{O}_{W} \longrightarrow 0,
\end{equation}

\noindent where $F \in \mathcal{R}(0,3,2s),\ W=\{x_1, ..., x_s\} \subset \mathbb{P}^3$. According to \cite{Chang} for a general sheaf $F$ in $\mathcal{R}(0,3,2s)$ there exists a nonzero section $\sigma \in H^{0}(F(1))$, the zero-scheme $Y=(\sigma)_0$ of which can be described as follows:

\begin{enumerate} [(i)]

\item for $s=1$, the scheme $Y$ is a disjoint union $l \sqcup l' \sqcup C_2$ of two lines $l, l'$ and a nonsingular conic $C_2$;

\item for $s=2$, the scheme $Y$ is a disjoint union of a line $l$ and a nonsingular twisted cubic $C_{3}$;

\item for $s=3$, the scheme $Y$ is a nonsingular rational quartic curve $C_{4}$;
\end{enumerate}
In all three cases (i)-(iii) the sheaf $F$ is a nontrivial extension of $\mathcal{O}_{\mathbb{P}^3}$-sheaves of the form
\begin{equation}\label{extn}
0 \longrightarrow \mathcal{O}_{\mathbb{P}^3}(-1) \longrightarrow F \longrightarrow \mathcal{I}_{Y}(1) \longrightarrow 0.
\end{equation}
Such extensions are classified by the vector space $V=Ext^{1}(\mathcal{I}_{Y}(1), \mathcal{O}_{\mathbb{P}^3}(-1)) \simeq H^0(\omega_{Y}(2))$,
and there exists a universal flat family of such extensions with base $V \simeq \mathbb{A}^{n}, n=\text{dim}~V$, - see
\cite[Proposition 3.1]{HH}. Now prove the following theorem.

\begin{Teo}\label{Thm2}
The following inclusions are true: $\overline{\mathcal{X}(1,1)} \subsetneqq \mathcal{T}(1), \ \overline{\mathcal{X}(2,2)} \subsetneqq \mathcal{T}(2),$ \\ $ \overline{\mathcal{X}(1,2)}  \subsetneqq \mathcal{T}(2), \ \overline{\mathcal{X}(2,3)} \subsetneqq \mathcal{T}(3).$
Hence the families $\overline{\mathcal{X}(1,1)}$,
$\overline{\mathcal{X}(2,2)}$, $\overline{\mathcal{X}(1,2)}$,
$\overline{\mathcal{X}(2,3)}$
do not constitute components in $\mathcal{M}(3)$.
\end{Teo}
\vspace{1mm}\noindent\textit{Proof.}
1) We first show that $\overline{\mathcal{X}(1,1)}\subsetneqq \mathcal{T}(1)$.
Here $\mathcal{T}(1)$ corresponds to the case (i) above,
where $Y=l \sqcup l' \sqcup C_2$. In this case $V \simeq H^0(\omega_{Y}(2))=  H^0(\mathcal{O}_{l}) \oplus H^0(\mathcal{O}_{l'}) \oplus H^0(\omega_{C_2}(2))\simeq \mathbb{A}^{5}$, and we can introduce the coordinates $(t_1,,...,t_5)$ on $\mathbb{A}^{5}$ such that $t_1\in H^0(\mathcal{O}_{l}),\ t_2\in H^0(\mathcal{O}_{l'}),\ (t_3,t_4,t_5)\in H^0(\omega_{C_2}(2))$. Let $\mathbf{F}=\{F_t\}_{t \in \mathbb{A}^1}$ be the flat subfamily of the universal family of extensions (\ref{extn})
over the affine line $\mathbb{A}^1=\{(t,1,...,1)\in \mathbb{A}^{5}~| ~ t \in \mathbf{k}\}$:
\begin{equation}\label{family bf F}
0 \longrightarrow \mathcal{O}_{\mathbb{P}^3}(-1) \boxtimes \mathcal{O}_{\mathbb{A}^1} \longrightarrow {\bf F} \longrightarrow \mathcal{I}_{{\bf Y} / \mathbb{P}^3 \times \mathbb{A}^1} \otimes \mathcal{O}_{\mathbb{P}^3}(1) \boxtimes \mathcal{O}_{\mathbb{A}^1} \longrightarrow 0,
\end{equation}
where $\mathbf{Y}=Y\times\mathbb{A}^1$. By construction,
\begin{equation}\label{in R(032)}
[F_t]\in\mathcal{R}(0,3,2),\ \ \ \ \ t\ne0,
\end{equation}
and, in addition, the reflexive sheaf $F_{0}^{\vee \vee}=
(\mathbf{F}|_{\mathbb{P}^3\times\{0\}})^{\vee \vee}$ fits into the short exact sequences
\begin{equation}\label{F0 dualdual}
0 \longrightarrow \mathcal{O}_{\mathbb{P}^3}(-1) \longrightarrow F_{0}^{\vee \vee} \longrightarrow \mathcal{I}_{Y'}(1) \longrightarrow 0,\ \ \ \
Y'=l'\sqcup C_2.
\end{equation}
\begin{equation}\label{F0}
0\longrightarrow F_{0} \longrightarrow F_{0}^{\vee \vee} \longrightarrow \mathcal{O}_l(1) \longrightarrow 0.
\end{equation}
Now (\ref{F0 dualdual}) yields:
\begin{equation}\label{in R(022)}
[F_{0}^{\vee \vee}]\in\mathcal{R}(0,2,2).
\end{equation}
Next, $W=\{q\}$ is a point in $\mathbb{P}^3$ such that, for each $t\in\mathbb{A}^1$, the sheaf $F_t$ of the family $\mathbf{F}$ is locally free at the point $W$. It follows that $\mathbf{F}|_{\mathbf{W}}\simeq 2\mathcal{O}_{\mathbf{W}}$,
where $\mathbf{W}=W\times\mathbb{A}^1$. Therefore, there
exists an epimorphism of sheaves $\mathbf{e}:\mathbf{F}\twoheadrightarrow
\mathcal{O}_{\mathbf{W}}$, and we obtain a flat family of sheaves
$\mathbf{E}=\{E_t\}_{t \in \mathbb{A}^1}$ defined by the exact triple:
\begin{equation}\label{family bf E}
0 \longrightarrow \mathbf{E}\longrightarrow {\bf F} \overset{\mathbf{e}}{\longrightarrow}\mathcal{O}_{\mathbf{W}} \longrightarrow 0.
\end{equation}
Restricting this triple onto $\mathbb{P}^3\times\{t\},\  t\in\mathbb{A}^1$, we obtain the triple
\begin{equation}\label{Et}
0\longrightarrow E_t \longrightarrow F_t \longrightarrow \mathcal{O}_{W} \longrightarrow 0,
\end{equation}
which together with
(\ref{E in T1}) and (\ref{in R(032)}) yields $[E_t]\in\mathcal{T}(1)$ for $t\ne0$. Hence also
\begin{equation}\label{E0 in T1}
[E_0]\in\mathcal{T}(1).
\end{equation}
Now (\ref{F0}) and the triple (\ref{Et}) for $t=0$ yield the exact triple
\begin{equation}\label{triple 2 for E0}
0 \longrightarrow E_0 \longrightarrow F_{0}^{\vee \vee} \longrightarrow \mathcal{O}_{l}(1) \oplus \mathcal{O}_{W} \longrightarrow 0,
\end{equation}
which together with (\ref{in R(022)}) implies that $[E_0]\in\mathcal{X}(1,1)$. Therefore, by (\ref{E0 in T1})
\begin{equation}\label{intersect}
[E_0]\in\mathcal{X}(1,1)\cap\mathcal{T}(1).
\end{equation}
Since
\begin{equation}\label{dims}
\dim\mathcal{T}(1)=25=\dim \mathrm{Ext}^1(E_0,E_0)=\dim  T_{[E_0]}\mathcal{M}(3),
\end{equation}
(see Section 1 and the above table)
it follows that $[E_0]$ is a smooth point of the moduli scheme $\mathcal{M}(3)$, we conclude the inclusion $\overline{\mathcal{X}(1,1)} \subset\mathcal{T}(1)$. As $\dim\mathcal{X}(1,1)=24<25,$ it follows that $\overline{\mathcal{X}(1,1)} \subsetneqq \mathcal{T}(1)$.

2) We next show that
$\overline{\mathcal{X}(2,2)}\subsetneqq \mathcal{T}(2)$.
Here $\mathcal{T}(2)$ corresponds to the case (ii) above,
where $Y=l \sqcup C_3$. In this case $V \simeq H^0(\omega_{Y}(2))=  H^0(\mathcal{O}_{l}) \oplus H^0(\omega_{C_3}(2))\simeq \mathbb{A}^6$, and we can introduce the coordinates $(t_1,,...,t_6)$ on $\mathbb{A}^{6}$ such that $t_1\in H^0(\mathcal{O}_{l}),\ (t_2,...,t_6)\in H^0(\omega_{C_3}(2))$. Let $\mathbf{F}=\{F_t\}_{t \in \mathbb{A}^1}$ be the flat subfamily of the universal family of extensions (\ref{extn})
over the affine line $\mathbb{A}^1=\{(t,1,...,1)\in \mathbb{A}^6~| ~ t \in \mathbf{k}\}$, so again $\mathbf{F}$
fits in the triple (\ref{family bf F}).
Furthermore, instead of (\ref{in R(032)}) and
(\ref{F0 dualdual}) one has
\begin{equation}\label{in R(034)}
[F_t]\in\mathcal{R}(0,3,4),\ \ \ \ \ t\ne0,
\end{equation}
\begin{equation}\label{F0 dualdual new}
0 \longrightarrow \mathcal{O}_{\mathbb{P}^3}(-1) \longrightarrow F_{0}^{\vee \vee} \longrightarrow \mathcal{I}_{C_3}(1) \longrightarrow 0,
\end{equation}
and the triple (\ref{F0}) remains true.
Now (\ref{F0 dualdual new}) yields:
\begin{equation}\label{in R(024)}
[F_{0}^{\vee \vee}]\in\mathcal{R}(0,2,4).
\end{equation}
Next, $W=\{q_1,q_2\}$ is now a disjoint union of two points in $\mathbb{P}^3$, and again the triples (\ref{family bf E}), (\ref{Et}) for $t\in U$, and (\ref{triple 2 for E0}) remain true.
Now using (\ref{in R(034)})-(\ref{in R(024)}) and arguing as above, we obtain similar to (\ref{intersect})and (\ref{dims}) that $[E_0]\in\mathcal{X}(2,2)\cap\mathcal{T}(2)$ and
$\dim\mathcal{T}(2)=29=\dim  T_{[E_0]}\mathcal{M}(3)$.
As $\dim\mathcal{X}(2,2)=26<29,$ it follows that $\overline{\mathcal{X}(2,2)} \subsetneqq \mathcal{T}(2)$.

3) Show that
$\overline{\mathcal{X}(1,2)}\subsetneqq \mathcal{T}(2)$.
For this, consider a 1-dimensional flat family $\pi:{\bf Y}\to U$ of curves in $\mathbb{P}^3$ with base $U$ with
marked point 0 being an open subset of $\mathbb{A}^1$, such that

\begin{enumerate} [(i)]

\item for $0\ne t\in U$ the fibre $Y_t$ of the family ${\bf Y}$ is a disjoint union $l \sqcup C_{3,t}$ of a line $l$ and a nonsingular twisted cubic $C_{3,t}$;

\item the zeroth fiber $Y_0$ of this family is a disjoint union $Y_0= l\sqcup Y'$, where $Y'$ is a union $Y'=l'\cup C_2$ of a line $l'$ and a nonsingular conic $C_2$, intersecting transversely at one point,
say, $p$. This yields an exact sequence
\end{enumerate}
\begin{equation}\label{omega(2)}
0\longrightarrow \omega_{C_2}(2)\oplus\omega_l(2)
\overset{i}{\longrightarrow}\omega_{Y_0}(2)\longrightarrow
\omega_{l'}(2) \oplus \mathcal{O}_{p}\longrightarrow0.
\end{equation}
On $U$ there is the vector bundle $\mathcal{V}=
\mathcal{E}xt^{1}_{\pi}(\mathcal{I}_{\mathbf{Y},\mathbb{P}^3\times U}\otimes\mathcal{O}_{\mathbf{P}^3}(1)\boxtimes
\mathcal{O}_U,\mathcal{O}_{\mathbb{P}^3}(-1)\boxtimes
\mathcal{O}_U)\simeq\pi_*(\omega_{\mathbf{Y},\mathbb{P}^3
\times U} \otimes \mathcal{O}_{\mathbb{P}^3}(2)\boxtimes
\mathcal{O}_U)\simeq6\mathcal{O}_U$.
Any section $s\in H^0(\mathcal{V})\simeq\mathrm{Ext}^{1}(\mathcal{I}_{\mathbf{Y},\mathbb{P}^3\times U}\otimes\mathcal{O}_{\mathbf{P}^3}(1)\boxtimes
\mathcal{O}_U,\mathcal{O}_{\mathbb{P}^3}(-1)\boxtimes
\mathcal{O}_U)$ defines a family $\mathbf{F}=\{F_t\}_{t \in U}$ of sheaves as an extension
\begin{equation}\label{family F 2}
0 \longrightarrow \mathcal{O}_{\mathbb{P}^3}(-1) \boxtimes \mathcal{O}_{U} \longrightarrow {\bf F} \longrightarrow \mathcal{I}_{{\bf Y} / \mathbb{P}^3 \times U} \otimes \mathcal{O}_{\mathbb{P}^3}(1) \boxtimes \mathcal{O}_{U} \longrightarrow 0.
\end{equation}
By flat base-change the section $s$ may be considered as a family of sections
$s=\{s_t\in H^0(\omega_{Y_t}(2))\}_{t\in U}$.
Now take a section $s$ such that
\begin{equation}\label{s0}
s_0|_{l'}=0,\ \ \ s_0|_{l}\ne0,\ \ \ s_0|_{C_2}\ne0,\ \ \ s_t|_{l}\ne0,\ \ \
s_t|_{C_{3,t}}\ne0,\ \ \ t\in U\setminus\{0\}
\end{equation}
(By (\ref{omega(2)}), to ensure that $s_0|_{l'}=0$
is the same as to take
$s_0\in H^0(i(\omega_{C_2}(2)\oplus\omega_l(2)))$).
Using the condition $s_0|_{l'}=0$, similarly to
(\ref{F0}) we obtain from (\ref{family F 2}) the triples
\begin{equation}\label{F0 2}
0\longrightarrow F_{0} \longrightarrow F_{0}^{\vee \vee} \longrightarrow \mathcal{O}_{l'}\longrightarrow 0,
\end{equation}
\begin{equation}\label{F0 dualdual 2}
0 \longrightarrow \mathcal{O}_{\mathbb{P}^3}(-1) \longrightarrow F_{0}^{\vee \vee} \longrightarrow \mathcal{I}_{l\sqcup C_2}(1) \longrightarrow 0.
\end{equation}
The triple (\ref{F0 dualdual 2}) yields (\ref{in R(022)}).
Besides, the exact triple (\ref{family bf E}) is true with
$\mathbf{W}=\{q_1,q_2\}\times U$. Arguing now as in the case 2), we obtain the inclusion $\overline{\mathcal{X}(1,2)}\subsetneqq \mathcal{T}(2)$.

4) Finally, show that $\overline{\mathcal{X}(2,3)} \subsetneqq \mathcal{T}(3)$. Consider a 1-dimensional flat family $\pi:{\bf Y}\to U$ of curves in $\mathbb{P}^3$ with base $U$ with marked point 0 being an open subset of $\mathbb{A}^1$, such that

\begin{enumerate} [(i)]

\item for $0\ne t\in U$ the fibre $Y_t$ of the family ${\bf Y}$ is a smooth rational quartic curve;

\item the zeroth fiber $Y_0$ of this family is a union $Y_0=l\cup C_3$ of a line $l$ and a nonsingular twisted cubic $C_3$ intersecting transversely at one point, say, $p$.
\end{enumerate}
As in case 2) above, on $U$ there is the vector bundle $\mathcal{V} \simeq \pi_*(\omega_{\mathbf{Y},\mathbb{P}^3
\times U} \otimes \mathcal{O}_{\mathbb{P}^3}(2)\boxtimes
\mathcal{O}_U)\simeq7\mathcal{O}_U$.
Again, a section $s\in H^0(\mathcal{V})$ defines a family $\mathbf{F}=\{F_t\}_{t \in U}$ of sheaves as an extension
(\ref{family F 2}), and this section $s$ may be considered as a family of sections $s=\{s_t\in H^0(\omega_{Y_t}(2))\}_{t\in U}$.
Taking a section $s$ such that
\begin{equation}\label{s0 new}
s_0|_{l}=0,\ \ \ s_0|_{C_3}\ne0,\ \ \ s_t|_{Y_t}\ne0,\ \ \ t\in U\setminus\{0\},
\end{equation}
we obtain similar to (\ref{F0 2}) and (\ref{F0 dualdual 2})
the exact triples
$0\longrightarrow F_{0} \longrightarrow F_{0}^{\vee \vee} \longrightarrow \mathcal{O}_l\longrightarrow 0$ and
$0 \longrightarrow \mathcal{O}_{\mathbb{P}^3}(-1) \longrightarrow F_{0}^{\vee \vee} \longrightarrow \mathcal{I}_{C_3}(1) \longrightarrow 0$,
Besides, the exact triple (\ref{family bf E}) is true with
$\mathbf{W}=\{q_1,q_2,q_3\}\times U$. Arguing now as in the case 3), we obtain the inclusion $\overline{\mathcal{X}(2,3)}\subsetneqq\mathcal{T}(3)$.
\hfill$\Box$

\end{document}